\documentclass[a4paper]{article}     
\usepackage{theorem}                 
\usepackage{amssymb}
\usepackage{epsfig}
\newtheorem{theorem}{Theorem}

\newcommand{\CC}{{\mathbb{C}}}

\newcommand{\RR}{{\mathbb{R}}}
\newcommand{\ZZ}{{\mathbb{Z}}}

\begin{document}

\title{Strange duality and polar duality}
\author{Wolfgang Ebeling}
\date{}

\maketitle

\begin{abstract}
We describe a relation between Arnold's strange duality and a
polar duality between the Newton polytopes which is mostly due to
M.~Koba\-ya\-shi. We show that this relation continues to hold for the
extension of Arnold's strange duality found by C.~T.~C.~Wall and the author.
By a method of Ehlers-Varchenko,  the characteristic polynomial of the
monodromy of a hypersurface singularity  can be computed from the Newton
diagram.  We generalize this method to the isolated
complete intersection singularities embraced in the extended duality. We use
this to explain the duality
of characteristic polynomials of the monodromy discovered by K.~Saito for
Arnold's original strange duality and extended by the author to the other
cases.
\end{abstract}

\section*{Introduction}

In \cite{Ebeling96} we gave a survey on some new features of Arnold's strange
duality and the extension of it found by C.~T.~C.~Wall and the author
\cite{EW85}. This duality can be considered as a
two-dimensional analogue of the Mirror Symmetry of Calabi-Yau threefolds.
Among other things, we discussed  a duality of the characteristic polynomials
of the monodromy operators of the singularities discovered by K.~Saito
\cite{Saito94}. We  showed that Saito's duality continues to hold for our
extension of Arnold's strange duality.

V.~Batyrev \cite{Batyrev94} showed that the Mirror Symmetry of Calabi-Yau
hypersurfaces in toric varieties is related to the polar duality between
their Newton polytopes. M.~Koba\-ya\-shi \cite{Kobayashi95} observed that
Arnold's strange duality corresponds to a duality of weight systems and this
in turn is related to a polar duality between certain polytopes.
These polytopes are slightly modified Newton polytopes of the singularities.

In this paper we
consider this relation more closely. We give a more precise and simpler
formulation of Kobayashi's result. We consider this correspondence for our
extension of Arnold's strange duality. This extension embraces isolated
complete intersection singularities (abbreviated ICIS in the sequel). We
define Newton diagrams and Newton polyhedra for these singularities and show
that Kobayashi's correspondence continues to hold for our extension of the
duality.

By a method found by A.~N.~Varchenko \cite{Varchenko76} and independently by
F.~Ehlers \cite{Ehlers78}, the characteristic polynomial of the
monodromy of a hypersurface singularity can be computed from the Newton
diagram. A generalization of this method to the principal monodromy of
non-degenerate complete intersection singularities has been given by M.~Oka
\cite{Oka90}. We show that for all the isolated complete intersection
singularities embraced in the duality with one exception, the
Ehlers-Varchenko method can be applied more directly. We use this to show
that the polar duality between the Newton polyhedra leads to Saito's duality
of characteristic polynomials.

\section{Newton diagrams}
\label{Sect1}

We start by fixing some basic notations.

For a set $M \subset \RR^n$, the {\em positive hull} of $M$ or the {\em cone}
determined by $M$ is the set of all non-negative linear combinations of
elements of $M$ and is denoted by $\mbox{pos}\, M$.

An $(n+1)$-tuple ${\bf w}:=(w_1,
\ldots ,w_n; N)$ of positive integers with $N \in \ZZ w_1 + \ldots + \ZZ w_n$
is called a {\em weight system}. The integers $w_i$ are called the {\em
weights} and $N$ is called the {\em degree} of ${\bf w}$.  A weight system
${\bf w}=(w_1, \ldots ,w_n; N)$ is called {\em reduced} if the weights
$w_1, \ldots , w_n$ have no common divisor. In the sequel we shall only
consider reduced weight systems.

Let
$$p(z)= \sum_{\nu \in \ZZ^n} a_\nu z^\nu \in \CC[z_1,z_1^{-1}, \ldots ,
z_n,z_n^{-1}]$$
be a Laurent polynomial. The {\em support} $\mbox{supp}(p)$ is the set
$$ \{ \nu \in \ZZ^n | a_\nu \neq 0 \}. $$
The Laurent polynomial $p$ is called {\em weighted homogeneous} if there exist
positive integers $w_1, \ldots , w_n$ (called {\em weights}) and a positive
integer $N$ (called {\em degree}) such that $ w_1 \nu_1 + \ldots + w_n \nu_n
= N$ for all $\nu \in \mbox{supp}(p)$. Then ${\bf w}:=(w_1, \ldots ,w_n; N)$
is a weight system.
Without loss of generality we may and will assume that $\bf w$ is reduced.
We call it the {\em weight system associated with} $p$.

The correspondence considered in \cite{Ebeling96} embraces the following
singularities:

\begin{itemize}
\item[a)] 14 exceptional unimodal hypersurface singularities given by
germs of analytic functions $f : (\CC^3,0) \to (\CC,0)$: $E_{12}$,
$E_{13}$, $E_{14}$, $Z_{11}$, $Z_{12}$, $Z_{13}$, $Q_{10}$, $Q_{11}$, $Q_{12}$,
$W_{12}$, $W_{13}$, $S_{11}$, $S_{12}$, $U_{12}$.
\item[b)] 6 bimodal quadrilateral hypersurface singularities given by
germs of analytic functions $f : (\CC^3,0) \to (\CC,0)$: $J_{3,0}$,
$Z_{1,0}$, $Q_{2,0}$, $W_{1,0}$, $S_{1,0}$, $U_{1,0}$.
\item[c)] 7 ICIS given by germs of analytic mappings $(g,f): (\CC^4,0) \to
(\CC^2,0)$, where $g$ and $f$ are given by
\begin{eqnarray*}  g(x,y,z,w) & = & xw-y^2 \\
                   f(x,y,z,w) & = & f'(x,y,w) + z^2
\end{eqnarray*}
for an analytic function $f': \CC^3 \to \CC$: $J'_9$,
$J'_{10}$,
$J'_{11}$, $K'_{10}$, $K'_{11}$, $J'_{2,0}$, $K'_{1,0}$.
\item[d)] 5 ICIS given by germs of analytic mappings $(g,f): (\CC^4,0) \to
(\CC^2,0)$, where $g$ is given by
\begin{eqnarray*}  g(x,y,z,w) & = & xw-yz:
\end{eqnarray*}
$L_{10}$, $L_{11}$, $M_{11}$, $L_{1,0}$, $M_{1,0}$.
\item[e)] the ICIS $I_{1,0}$ given by the equations
\begin{eqnarray*} g(x,y,z,w) & = & x^3 - yw \\
                  f(x,y,z,w) & = & (a+1)x^3 + yz +z^2 + zw
\end{eqnarray*}
where $a \neq 0,1$.
\end{itemize}

The precise equations can be found in \cite{Ebeling96} (see also below).

We shall now define Newton diagrams for these singularities. We recall the
definition of the Newton diagram of a Laurent polynomial.

Let  $p \in \CC[z_1,z_1^{-1}, \ldots , z_n,z_n^{-1}]$ be a Laurent
polynomial.

\addvspace{3mm}

\noindent {\bf Definition }(cf.\ \cite{Kouchnirenko76}) The convex hull in
$\RR^n$ of the set $\mbox{supp}(p) \setminus \{ 0 \}$ is called the {\em
Newton diagram} $\Gamma(p)$ of $p$.

\addvspace{3mm}

Let $p$ be weighted homogeneous with weight system ${\bf w}:=(w_1, \ldots
,w_n; N)$. The Newton diagram of $p$ is an
$(n-1)$-dimensional convex polytope. It can also be characterized in the
following way. Let $\RR_+$ be the set
of all non-negative real numbers. We denote by $\Gamma_+(p)$ the convex hull of
the set
$$ \bigcup_{\nu \in \, {\scriptstyle \rm supp}(p)} (\nu + \RR_+^n).$$
Then $\Gamma(p)$ is the union of all compact faces of $\Gamma_+(p)$. The
Newton diagram $\Gamma(p)$ lies in the hyperplane
$$ w_1 \nu_1 + \ldots + w_n \nu_n = N.$$

In order to define Newton diagrams for the above singularities,
we associate a Laurent polynomial in three variables
$$p(x,y,z)=\sum_{\nu \in \ZZ^3} a_\nu x^{\nu_1}y^{\nu_2}z^{\nu_3} \in \CC [
x,x^{-1},y,y^{-1},z,z^{-1}] $$
to each of these singularities. In each case we can assume that the function
$f$ is a weighted homogeneous polynomial with weight system ${\bf w}=(w_1,
\ldots , w_n;N)$ where $n=3$ or $n=4$. In the cases a) and b) we simply take
$p:=f$. In the cases c) and d)  we make the following substitution in
$f(x,y,z,w)$ to obtain a weighted homogeneous Laurent polynomial $p(x,y,z)$:
\begin{itemize}
\item[c)] $w=x^{-1}y^2$
\item[d)] $w=x^{-1}yz$
\end{itemize}
These equations are simply obtained by eliminating the variable $w$ from the
corresponding equations $g=0$.

The case e) is somewhat exceptional since it is the only case among the ICIS
where the singularity cannot be given by a germ of an analytic mapping
$(g,f): (\CC^4,0) \to (\CC^2,0)$ with first equation defining an
$A_1$-singularity. Moreover, the function $g$ in e) has a non-isolated
singularity in the origin, and the singularity $I_{1,0}$ is not a stabilization
of a curve singularity as in c). We might as well substitute $w=x^3y^{-1}$ in
the polynomial $f(x,y,z,w)$ of e), but unfortunately, it turns out that this
doesn't yield the same results as in the other cases. Instead of that  we
associate the following Laurent polynomial to this singularity:
$$p(x,y,z):=x^3+x^3y^{-3}z^3+yz.$$
Two terms also occur in the polynomial $f(x,y,z,w)$ of e) as above, but a
third term $x^3y^{-3}z^3$ appears which we cannot explain.

The weight systems
and Laurent polynomials are listed in Table~\ref{Table1}.

\begin{table} \centering
\caption{Weight systems and Laurent polynomials} \label{Table1}
\begin{tabular}{|c|c|c|c|}  \hline
Name  & ${\bf w}$ & $p$   & Dual
\\ \hline

$E_{12}$ & $(6,14,21;42)$ & $x^7+y^3+z^2$ & $E_{12}$  \\
$E_{13}$ & $(4,10,15;30)$ & $x^5y+y^3+z^2$ & $Z_{11}$ \\
$E_{14}$ & $(3,8,12;24)$ & $x^4z+y^3+z^2$ & $Q_{10}$ \\
$Z_{11}$ & $(6,8,15;30)$ & $x^5+xy^3+z^2$ & $E_{13}$ \\
$Z_{12}$ & $(4,6,11;22)$ & $x^4y+xy^3+z^2$ & $Z_{12}$ \\
$Z_{13}$ & $(3,5,9;18)$ & $x^3z+xy^3+z^2$ & $Q_{11}$ \\
$Q_{10}$ & $(6,8,9;24)$ & $x^4+y^3+xz^2$ & $E_{14}$ \\
$Q_{11}$ & $(4,6,7;18)$ & $x^3y+y^3+xz^2$ & $Z_{13}$ \\
$Q_{12}$ & $(3,5,6;15)$ & $x^3z+y^3+xz^2$ & $Q_{12}$ \\
$W_{12}$ & $(4,5,10;20)$ & $x^5+z^2+y^2z$ & $W_{12}$ \\
$W_{13}$ & $(3,4,8;16)$ & $x^4y+z^2+y^2z$ & $S_{11}$ \\
$S_{11}$ & $(4,5,6;16)$ & $x^4+xz^2+y^2z$ & $W_{13}$ \\
$S_{12}$ & $(3,4,5;13)$ & $x^3y+xz^2+y^2z$ & $S_{12}$ \\
$U_{12}$ & $(3,4,4;12)$ & $x^4+y^2z+yz^2$ & $U_{12}$ \\
\hline
$J_{3,0}$ & $(2,6,9;18)$ & $x^6y+y^3+z^2$ & $J'_9$\\
$Z_{1,0}$ & $(2,4,7;14)$ & $x^5y+xy^3+z^2$ & $J'_{10}$ \\
$Q_{2,0}$ & $(2,4,5;12)$ & $x^4y+y^3+xz^2$ & $J'_{11}$ \\
$W_{1,0}$ & $(2,3,6;12)$ & $x^6+z^2+y^2z$ & $K'_{10}$ \\
$W_{1,0}$ & $(2,3,6;12)$ & $x^3y^2+x^3z+y^2z+z^2$ & $L_{10}$ \\
$S_{1,0}$ & $(2,3,4;10)$ & $x^5+xz^2+y^2z$ & $K'_{11}$ \\
$S_{1,0}$ & $(2,3,4;10)$ & $x^2y^2+x^3z+y^2z+xz^2$ & $L_{11}$ \\
$U_{1,0}$ & $(2,3,3;9)$ & $x^3z+x^3y+yz^2+y^2z$ & $M_{11}$ \\
\hline
$J'_9$ & $(6,8,9,10;18)$ & $x^3+x^{-1}y^3+z^2$ & $J_{3,0}$ \\
$J'_{10}$ & $(4,6,7,8;14)$ & $x^2y+x^{-1}y^3+z^2$ & $Z_{1,0}$ \\
$J'_{11}$ & $(3,5,6,7;12)$ & $x^2z+x^{-1}y^3+z^2$ & $Q_{2,0}$ \\
$K'_{10}$ & $(4,5,6,6;12)$ & $x^3+z^2+x^{-1}y^2z$ & $W_{1,0}$ \\
$K'_{11}$ & $(3,4,5,5;10)$ & $x^2y+z^2+x^{-1}y^2z$ & $S_{1,0}$ \\
$L_{10}$ & $(4,5,6,7;12)$ & $x^3+z^2+x^{-1}y^2z$ & $W_{1,0}$ \\
$L_{11}$ & $(3,4,5,6;10)$ & $x^2y+z^2+x^{-1}y^2z$ & $S_{1,0}$ \\
$M_{11}$ & $(3,4,4,5;9)$ & $x^3+x^{-1}y^2z+x^{-1}yz^2$ & $U_{1,0}$ \\
\hline
$J'_{2,0}$ & $(2,4,5,6;10)$ & $x^3y+x^{-1}y^3+z^2$ & $J'_{2,0}$ \\
$L_{1,0}$ & $(2,3,4,5;8)$ & $xy^2+x^2z+x^{-1}y^2z+z^2$ & $L_{1,0}$ \\
$L_{1,0}$ & $(2,3,4,5;8)$ & $x^4+z^2+x^{-1}y^2z$ & $K'_{1,0}$ \\
$K'_{1,0}$ & $(2,3,4,4;8)$ & $xy^2+x^2z+x^{-1}y^2z+z^2$ & $L_{1,0}$ \\
$K'_{1,0}$ & $(2,3,4,4;8)$ & $x^4+z^2+x^{-1}y^2z$ & $K'_{1,0}$ \\
$M_{1,0}$ & $(2,3,3,4;7)$ & $x^2z+x^2y+x^{-1}yz^2+x^{-1}y^2z$ & $M_{1,0}$ \\
$I_{1,0}$ & $(2,3,3,3;6)$ & $x^3+x^3y^{-3}z^3+yz$ & $I_{1,0}$ \\
\hline

\end{tabular}
\end{table}

M.~Kobayashi \cite{Kobayashi95} has defined a polytope associated to a weight
system. We modify his definition slightly.

Let ${\bf w}=(w_1, \ldots ,w_n;
N)$ be a reduced weight system. The hyperplane $w_1 \nu_1 +
\ldots + w_n \nu_n = N$ meets the coordinate axes in the points $N/w_1, \ldots
, N/w_n$. We consider the convex hull of these $n$ points.

\addvspace{3mm}

\noindent {\bf Definition } The {\em full Newton
diagram} $\Delta({\bf w})$ of ${\bf w}$ is the $(n-1)$-simplex in $\RR^n$ which
is the convex hull of the vertices $(\frac{N}{w_1},0, \ldots , 0)$ ,
$\ldots$ , $(0, \ldots , 0, \frac{N}{w_n})$.

\addvspace{3mm}

For the weight systems of the singularities of classes a) and b) we have $n=3$.
In this case the full Newton diagram is a triangle.  In the
remaining cases we have weight systems ${\bf w}$ with $n=4$. We
shall now define for the weight system of a singularity of class c), d), or e)
a full Newton diagram $\Delta({\bf w})$.

First consider the weight system ${\bf w}=(w_1,w_2,w_3,w_4;N)$ of a singularity
of class c). Here
$g(x,y,z,w)=xw-y^2$. This implies that the weights $w_1$, $w_2$, $w_3$, $w_4$
have to satisfy the equation $w_1 + w_4 = 2w_2$. The substitution
$w=x^{-1}y^2$ means that we get an extra vector $(-1,2,0)^t$ in $\RR^3$ with
the coordinates $x$, $y$, $z$. Let $e_1$, $e_2$, $e_3$ be the vectors of the
standard basis of $\RR^3$.
Let
$\check{\sigma} := \mbox{pos}\{ e_1,-e_1+2e_2 \} \subset \RR^2$. According to
\cite[VI, 8, Example 1]{Ewald96} the affine toric variety
$X_{\check{\sigma}}$ corresponding to $\check{\sigma}$  is the quadric
$xw=y^2$.

\addvspace{3mm}

\noindent {\bf Definition }  The {\em full Newton diagram}
$\Delta({\bf w})$ of ${\bf w}$ is the triangle in $\RR^3$ which is the convex
hull of the  vertices
$(\frac{N}{w_1} , 0 , 0)$,
$(0, 0 , \frac{N}{w_3})$, and
$( - \frac{N}{w_4} , 2\frac{N}{w_4}, 0)$.

\addvspace{3mm}

Note that the point $(0,\frac{N}{w_2}, 0)$ lies in the
convex hull of the points $(\frac{N}{w_1}, 0, 0)$ and $(-
\frac{N}{w_4}, 2\frac{N}{w_4}, 0)$ because of the relation $w_1 + w_4 =
2w_2$.

Now consider the weight system ${\bf w}=(w_1,w_2,w_3,w_4;N)$ of a singularity
of class d). Here $g(x,y,z,w)=xw-yz$. From this we get the relation $w_1 + w_4
= w_2 + w_3$. The substitution $w=x^{-1}yz$ means that we get an extra vector
$(-1,1,1)^t$. Let $\check{\sigma} :=
\mbox{pos}\{ e_1,e_2,e_3,-e_1+e_2+e_3 \}$. As in \cite[loc.~cit.]{Ewald96}
one can see that the affine toric variety
$X_{\check{\sigma}}$ corresponding to $\check{\sigma}$  is the quadric
$xw=yz$.

\addvspace{3mm}

\noindent {\bf Definition }  The {\em full Newton diagram}
$\Delta({\bf w})$ of ${\bf w}$ is the quadrilateral in $\RR^3$ which is the
convex hull of the  points
$(\frac{N}{w_1} , 0 , 0)$, $(0,\frac{N}{w_2},0)$,
$(0, 0 , \frac{N}{w_3})$, and
$(-\frac{N}{w_4},\frac{N}{w_4},\frac{N}{w_4})$.

\addvspace{3mm}

Finally we consider the weight system ${\bf w} =
(w_1,w_2,w_3,w_4;N)=(2,3,3,3;6)$ of a singularity of type
$I_{1,0}$ (case e)). In this case we have $g(x,y,z,w)=x^3-yw$ and the relation
$3w_1=w_2+w_4$. Therefore we get an extra vector
$(3,-1,0)^t$. Let $\check{\sigma} = \mbox{pos} \{3e_1-e_2,e_2,e_3 \} \subset
\RR^3$. As in case c) one can see that the affine toric variety
$X_{\check{\sigma}}$ corresponding to $\check{\sigma}$  is given by the
equation $x^3=yw$.

\addvspace{3mm}

\noindent {\bf Definition }  The {\em full Newton diagram}
$\Delta({\bf w})$ of ${\bf w}$ is the triangle in $\RR^3$ which is the convex
hull of the vertices
$(3\frac{N}{w_4} , - \frac{N}{w_4} , 0)$,
$( 0, \frac{N}{w_2} , 0 )$, and
$( 0  , 0 , \frac{N}{w_3})$.

\addvspace{3mm}

Again note that the point $(\frac{N}{w_1} , 0, 0)$ lies in the convex
hull of the points $(3\frac{N}{w_4} , - \frac{N}{w_4} , 0)$ and
$( 0 , \frac{N}{w_2} , 0 )$.

\section{Polar duality}
\label{Sect2}

Now let $(X,0)$ be one of the above singularities, let
$p$ be the Laurent polynomial of $(X,0)$ according to
Table~\ref{Table1}, and let ${\bf w}$ be the weight system of $(X,0)$. The
Newton diagram of $(X,0)$ is
$\Gamma(p)$, and its full Newton diagram is
$\Delta({\bf w})$.

\addvspace{3mm}

\noindent {\bf Definition } (cf.\ \cite{Kobayashi95}) We shall call {\em
Newton polyhedron} (resp.\ {\em full Newton polyhedron}) and denote
by $\tilde{\Gamma}(p)$ (resp.\ $\tilde{\Delta}({\bf w})$) the polyhedron in
$\RR^3$ which is obtained from $\Gamma(p)$ (resp.\ $\Delta({\bf w})$)
by taking the convex hull with the origin and translating by the vector
$u_0:=(-1,-1,-1)$ in cases a) and b), $u_0:=(0,-1,-1)$ in cases c) and d),
and $u_0:=(-1,0,-1)$ in case e).

\addvspace{3mm}

The polyhedra $\tilde{\Gamma}:=\tilde{\Gamma}(p)$ and
$\tilde{\Delta}:=\tilde{\Delta}({\bf w})$ are both integral polyhedra
containing the
origin in their interiors. Let $u_0$, $u_1$, $u_2$, $u_3$ (,$u_4$) be the
vertices of $\tilde{\Gamma}$ with $u_0$ as above.

We shall now consider the polar duality between
polytopes.

\addvspace{3mm}

\noindent {\bf Definition } Let $M \subset \RR^n$. Let $\langle \ , \
\rangle$ denote the Euclidean scalar product of $\RR^n$. The {\em polar
dual} of $M$ is the following subset of $\RR^n$:
$$M^\ast := \{ y \in \RR^n | \langle x , y \rangle \geq -1 \ \mbox{for all}
\ x \in M \}.$$

\addvspace{3mm}

The polar dual $\tilde{\Delta}^\ast$ of $\tilde{\Delta}$ is again an integral
polyhedron. The vertices of
$\tilde{\Delta}^\ast$ can be computed in the different cases as follows:
\begin{itemize}
\item[a)] $v_1:=(1,0,0)$, $v_2:=(0,1,0)$, $v_3:=(0,0,1)$,
$v_0:=(-w_1,-w_2,-w_3)$,
\item[b)] same as a),
\item[c)] $v_1:=(2,1,0)$, $v_2:=(0,1,0)$, $v_3:=(0,0,1)$,
$v_0:=(-w_1,-w_2,-w_3)$,
\item[d)] $v_1:=(1,0,1)$, $v_2:=(1,1,0)$, $v_3:=(0,0,1)$, $v_4:=(0,1,0)$, \\
$v_0:=(-w_1,-w_2,-w_3)$.
\item[e)] $v_1:=(1,0,0)$, $v_2:=(1,3,0)$, $v_3:=(0,0,1)$,
$v_0:=(-w_1,-w_2,-w_3)$.
\end{itemize}

We translate the polytope
$\tilde{\Delta}^\ast$ by the vector
$v_0$ such that the new vertices are $\tilde{v}_i := v_i - v_0$, $i=1,2,3(,4)$
and the origin. The convex hull of the points $\tilde{v}_i$, $i=1,2,3(,4)$
will be denoted by $\nabla$.

An integral $n \times n$-matrix is called
{\em unimodular} if it has determinant $\pm 1$.

Our
main result is the following theorem.

\begin{theorem} \label{thm:main}
For each of the above singularities, there exists a
unique reduced weight system ${\bf w}^\ast$ with $N^\ast=N$ and a unique
unimodular
$3
\times 3$-matrix $A$ with entries in the positive integers which transforms
the polyhedron $\tilde{\Gamma}$ to the polyhedron $\tilde{\Delta}({\bf
w}^\ast)^\ast$ such that
$\{ u_1, u_2, u_3 \}$ (resp.\ $\{ u_1, u_2, u_3, u_4 \}$)
is mapped to $\{ v_1, v_2, v_3\}$ (resp.\ $\{ v_1, v_2, v_3, v_4 \}$). The
weight system ${\bf w}^\ast$ is the weight system of the dual singularity. In
the coordinate system given by taking the rows of $A$ as basis vectors of
$\RR^3$ the polygon
$\nabla$  is the Newton diagram of the dual singularity.
\end{theorem}

\noindent {\em Proof. }
Let $U$ be the matrix with column vectors $u_i$, and let $V$ be the matrix
with column vectors $v_i$, $i=1, \ldots , m$, $m=3$ or $m=4$.
Then either $m=4$ or $m=3$ and $|\det U| = 1,2,3$. The same holds for the
matrix $V$. Let $V_\Gamma$ be the matrix of the same size and if $m=3$, with
the same absolute value of the determinant as $U$. Then the matrix $A$ is
determined by the equation
$$ A U = V_\Gamma . $$
The weight system ${\bf w}^\ast$ is determined by $m$, $|\det U|$, $N$,
and the vector
$$-Au_0 = (w_1^\ast, w_2^\ast, w_3^\ast).$$
Namely,
$${\bf w}^\ast= \left\{ \begin{array}{ll}(w_1^\ast, w_2^\ast, w_3^\ast; N) &
\mbox{if } m=3 \mbox{ and } |\det U| =1, \\
(w_1^\ast,w_2^\ast, w_3^\ast, 2w_2^\ast-w_1^\ast; N) & \mbox{if } m=3 \mbox{
and }  |\det U| = 2, \\
(w_1^\ast,w_2^\ast, w_3^\ast, 3w_1^\ast-w_2^\ast; N) & \mbox{if } m=3 \mbox{
and }  |\det U| = 3, \\
(w_1^\ast,w_2^\ast, w_3^\ast, w_2^\ast+w_3^\ast-w_1^\ast; N) & \mbox{if } m=4.
\end{array} \right.$$
One checks case by case that
$A$ has entries in the positive integers and that
${\bf w}^\ast$ is the weight system and
$\nabla$ (in the coordinate system obtained by taking the rows of $A$ as basis
vectors) the Newton diagram of the dual singularity.
This proves Theorem~\ref{thm:main}.

\addvspace{3mm}

The matrices $A$ are listed in Table~\ref{Table2}. This table is to be read
as follows: Starting with a singularity to the left of $A$ one has to apply
$A$, starting with a singularity to the right of $A$ one has to apply $A^t$.

\begin{table} \centering
\caption{The matrices A} \label{Table2}
\begin{tabular}{|c|c|c||c|c|c|}  \hline
Name   & $A$    & Dual & Name & $A$ & Dual
\\ \hline

$E_{12}$  &
$\left( \begin{array}{ccc} 1 & 2 & 3 \\ 2 & 5 & 7 \\ 3 & 7 & 11
\end{array} \right)$  &
$E_{12}$
&

$E_{13}$ &
$\left( \begin{array}{ccc} 1 & 2 & 3 \\
1 & 3 & 4
\\ 2 & 5 & 8
\end{array} \right)$  &
$Z_{11}$
\\ \hline

$\begin{array}{c} E_{14} \\  J_{3,0} \end{array}$  &
$\left( \begin{array}{ccc} 1 & 2 & 3 \\ 1 & 3 & 4 \\ 1 & 3 & 5
\end{array} \right)$ &
$\begin{array}{c} Q_{10} \\  J'_9 \\  \end{array}$
&

$Z_{12}$ &
$\left( \begin{array}{ccc} 1 & 1 & 2 \\ 1 & 2 & 3 \\ 2 & 3 & 6
\end{array} \right)$ &
$Z_{12}$
\\ \hline

$\begin{array}{c} Z_{13} \\  Z_{1,0}  \end{array}$ &
$\left( \begin{array}{ccc} 1 & 1 & 2 \\ 1 & 2 & 3 \\ 1 & 2 & 4
\end{array} \right)$
& $\begin{array}{c}  Q_{11} \\ J'_{10}  \end{array}$
&

$\begin{array}{c} Q_{12} \\ Q_{2,0}  \\ J'_{2,0} \end{array}$ &
$\left( \begin{array}{ccc} 1 & 1 & 1 \\ 1 & 2 & 2 \\ 1 & 2 & 3
\end{array} \right)$   &
$\begin{array}{c} Q_{12} \\  J'_{11}  \\ J'_{2,0} \end{array}$
\\ \hline

$W_{12}$ &
$\left( \begin{array}{ccc} 1 & 1 & 2 \\ 1 & 1 & 3 \\ 2 & 3 & 5
\end{array} \right)$ &
$W_{12}$
&

$\begin{array}{c} W_{13} \\ W_{1,0} \\ W_{1,0}  \end{array}$  &
$\left( \begin{array}{ccc} 1 & 1 & 2 \\ 1 & 1 & 3 \\ 1 & 2 & 3
\end{array} \right)$ &
$\begin{array}{c}  S_{11} \\ K'_{10} \\ L_{10}  \end{array}$
\\ \hline

$\begin{array}{c} S_{12} \\ S_{1,0} \\ S_{1,0} \\
L_{1,0} \\ L_{1,0} \\ K'_{1,0} \\ I_{1,0} \end{array}$ &
$\left( \begin{array}{ccc} 1 & 1 & 1 \\ 1 & 1 & 2 \\ 1 & 2 & 2
\end{array} \right)$ &
$\begin{array}{c} S_{12} \\ K'_{11} \\ L_{11} \\
L_{1,0} \\ K'_{1,0} \\ K'_{1,0} \\ I_{1,0} \end{array}$
&

$\begin{array}{c} U_{12} \\ U_{1,0} \\ M_{1,0} \end{array}$ &
$\left( \begin{array}{ccc} 1 & 1 & 1 \\ 1 & 2 & 1 \\ 1 & 1 & 2
\end{array} \right)$ &
$\begin{array}{c} U_{12} \\ M_{11} \\  M_{1,0} \end{array}$
\\ \hline

\end{tabular}
\end{table}

\section{The characteristic polynomial of the mono\-dro\-my}
\label{Sect3}

Let $f: (\CC^n,0) \to (\CC,0)$ be a germ of an analytic function defining an
isolated hypersurface singularity $(X,0)$. Let $X_t$ denote the Milnor
fibre of $f$ over a point $t \in \CC$ with $|t| << 1$. A characteristic
homeomorphism of the Milnor fibration of $f$ induces an automorphism
$c:H_{n-1}(X_t,\ZZ) \to H_{n-1}(X_t,\ZZ)$ called the {\em (classical)
monodromy operator} of $(X,0)$. Let $\phi(\lambda) = \det (\lambda I - c)$ be
the characteristic polynomial of $c$.

By the Ehlers-Varchenko method \cite{Ehlers78, Varchenko76}, the polynomial
$\phi(\lambda)$ can be computed from the Newton diagram $\Gamma(f)$ as
follows:

Let
$$p(z) = \sum_{\nu \in \ZZ^n} a_\nu z^\nu \in \CC[z_1,z_1^{-1}, \ldots ,
z_n, z_n^{-1}]$$
be a Laurent polynomial. Let $s$ be a compact face of
$\Gamma(p)$. We define $p_s$ by $p_s(z):=\sum_{\nu \in s} a_\nu z^\nu$. We
call $p$ {\em non-degenerate} if for any compact
face $s \in \Gamma(p)$ the Laurent polynomials
$$z_1 \frac{\partial p_s}{\partial z_1}, \ldots , z_n \frac{\partial
p_s}{\partial z_n}$$
have no common zero in $\{ z \in \CC^n | z_1 \cdot \ldots
\cdot z_n \neq 0 \}$ (cf.\ \cite[1.19]{Kouchnirenko76}).

Let $b_1, \ldots , b_k$, $k \geq n$, be a system of
generators of $\RR^n$. We assume that $a_\nu = 0$ for $\nu \not\in
\mbox{pos}\{ b_1, \ldots , b_k \}$. For $I \subset \{1, \ldots , k \}$ let
$U_I$ be the linear span of the set $\{ b_i | i \in I \}$.
If an $r$-dimensional compact face $s$ of $\Gamma(p)$ lies entirely in an
$(r+1)$-dimensional  subspace $U_I$ for some index set $I$, it is
called {\em special}. For a special $r$-dimensional face
$s \subset U_I$ we define numbers $m(s)$, $V(s)$, and $\chi(s)$ as follows:
Let $\sum_{i \in I} \alpha_i x_i = m(s)$ be the equation of $s$ in $U_I$ with
respect to a basis of $U_I$ in $\{ b_i | i \in I \}$ where the $\alpha_i$ and
$m(s)$ are non-negative integers and the greatest common divisor of the numbers
$\alpha_i$, $i \in I$, is equal to one. Let
$V(s)$ be the $(r+1)$-dimensional volume of the convex hull of $\{ 0 \} \cup s$
multiplied by $(r+1)!$. We set
$$\chi(s) := (-1)^r \frac{V(s)}{m(s)}.$$

\addvspace{3mm}

\noindent {\bf Definition }
$$\phi_{\Gamma(p)}(\lambda) := \left( (\lambda -1)^{-1} \prod_{s \
{\scriptstyle \rm special}} (\lambda^{m(s)} -1)^{\chi(s)}
\right)^{(-1)^{n-1}}.$$

\addvspace{3mm}

\noindent Note that this definition depends on the choice of the system of
generators
$b_1, \ldots , b_k$ of $\RR^n$.

Let $f$ be given by a non-degenerate polynomial. By
\cite{Ehlers78, Varchenko76} we then have $\phi(\lambda)=
\phi_{\Gamma(f)}(\lambda)$  for $k=n$ and $b_1, \ldots , b_k$ being the
standard basis of $\RR^n$.

Now let $(X,0)$ be an ICIS defined by a germ of an analytic mapping $(g,f) :
(\CC^4,0) \to (\CC,0)$. We assume that $(g,f)$ are generically chosen
such that $(X',0):=(g^{-1}(0),0)$ is an
isolated hypersurface singularity of minimal Milnor number $\mu_1$ among such
choices of $g$. Then the {\em monodromy operator} of $(X,0)$ is defined to be
the monodromy operator of the function germ $f: (X',0) \to (\CC,0)$. We shall
show that the Ehlers-Varchenko method can be extended to the ICIS of the
classes c) and d).

First consider an ICIS $(X,0)$ given by
\begin{eqnarray*}  g(x,y,z,w) & = & xw-y^2 \\
                   f(x,y,z,w) & = & f'(x,y,w) + z^2.
\end{eqnarray*}
Define $g'(x,y,z) := xw-y^2$ and consider the curve singularity $(Y,0)$
defined by $(g',f') : (\CC^3,0) \to (\CC^2,0)$. Then
$(Y',0)=((g')^{-1}(0),0)$ is an isolated singularity of minimal Milnor number
$\mu_1 = 1$. Let $\phi'(\lambda)$ be the characteristic polynomial of the
monodromy operator $c'$ of $(Y,0)$. As above, we make the substitution
$w=x^{-1}y^2$ to obtain a Laurent polynomial
$p'(x,y) := f'(x,y,x^{-1}y^2)$. We assume that $p'$ is non-degenerate. Then
$\phi'(\lambda) =
\phi_{\Gamma(p')}(\lambda)$ where $\phi_{\Gamma(p')}(\lambda)$ is defined using
the basis $\{ b_1, b_2 \}$ where $b_1:=e_1$ and $b_2:=-e_1+2e_2$. This follows
from the fact that the affine toric variety
$X_{\check{\sigma}}$ corresponding to $\check{\sigma}:= \mbox{pos}\{
e_1,-e_1+2e_2 \}$  is the quadric
$xw=y^2$. Then the proof of Ehlers-Varchenko goes through replacing the
initial cone $\RR^n_+$ corresponding to $\CC^n$ by $\check{\sigma}$.

The original ICIS $(X,0)$ is a stabilization (2-fold suspension) of the ICIS
$(Y,0)$. By \cite{ESt96} the characteristic polynomial $\phi(\lambda)$ of the
monodromy of $(X,0)$ is related to the polynomial $\phi'(\lambda)$ by
replacing $\lambda$ by $-\lambda$ in the formula for $\phi'(\lambda)$ and
adding the term $(\lambda^2-1)$. This means that if
$$\phi'(\lambda) = (\lambda-1) \prod_{m > 1}(\lambda^m -1)^{\chi_m}$$
then
$$\phi(\lambda) = (\lambda^2-1)^2 (\lambda-1)^{-1} \prod_{m > 1}(\lambda^{2m}
-1)^{\chi_m} (\lambda^m -1)^{-\chi_m}.$$
This result can also be formulated in the following way:

\begin{theorem} \label{thm:EVc}
Let $(X,0)$ be an ICIS given by a germ of an analytic mapping $(g,f):
(\CC^4,0) \to (\CC^2,0)$ with
$g(x,y,z,w) =  xw-y^2$ and $f(x,y,z,w)  =  f'(x,y,w) + z^2$,
and let $p'(x,y):=f'(x,y,x^{-1}y^2)$, $p(x,y,z):=p'(x,y) + z^2$. Assume that
$p'$ is non-degenerate. Let $\phi_{\Gamma(p)}(\lambda)$ be defined using the
basis $\{e_1,-e_1+2e_2,e_3 \}$ of $\RR^3$. Then
$$\phi(\lambda) = (\lambda^2-1) \phi_{\Gamma(p)}(\lambda).$$
\end{theorem}

Similarly we can derive the following theorem.

\begin{theorem} \label{thm:EVd}
Let $(X,0)$ be an ICIS given by a germ of an analytic mapping $(g,f):
(\CC^4,0) \to (\CC^2,0)$ with
$g(x,y,z,w) =  xw-yz$, and let $p(x,y,z):=f(x,y,z,x^{-1}yz)$.
Assume that $p$ is
non-degenerate. Let
$\phi_{\Gamma(p)}(\lambda)$ be defined using the system of generators
$\{ e_1, e_2, e_3, -e_1+e_2+e_3 \}$ of $\RR^3$. Then
$$\phi(\lambda) = \phi_{\Gamma(p)}(\lambda).$$
\end{theorem}

\noindent {\em Proof. } This follows by an analogous argument as above
from the fact that the affine toric variety
$X_{\check{\sigma}}$ corresponding to $\check{\sigma}:=
\mbox{pos}\{ e_1,e_2,e_3,-e_1+e_2+e_3 \}$  is the quadric
$xw=yz$.

Finally consider the singularity $I_{1,0}$. This singularity is of the type
of \cite[Example (7.3)]{Oka90}. It follows from \cite[loc.~cit.]{Oka90} that
$$\phi(\lambda) = (\lambda -1)^{-1}(\lambda^2 -1)^{-2}(\lambda^6 -1)^3.$$

\section{Saito's duality of characteristic polynomials}
\label{Sect4}

Let $h$ be a positive integer and let $\psi(\lambda)$ be a polynomial in the
variable $\lambda$ of the form
$$\psi(\lambda)= \prod_{m | h} (\lambda^m -1)^{\chi_m} \quad \mbox{for} \
\chi_m \in \ZZ.$$
K.~Saito (\cite{Saito94}, cf.\ also \cite{Ebeling96}) has defined a {\em dual
polynomial}
$\psi^\ast(\lambda)$ to $\psi(\lambda)$:
$$\psi^\ast(\lambda) = \prod_{k | h} (\lambda^k -1)^{-\chi_{h/k}}.$$

Now let $(X,0)$ be a singularity of class a), b), c), or d). If $(X,0)$ is of
type a) or b), then we set
$$\psi(\lambda):=\phi(\lambda),$$
the characteristic polynomial of the monodromy. For $(X,0)$ of type c) or d)
we put
$$\psi(\lambda):=(\lambda -1)^{-1}\phi(\lambda).$$
By Sect.~\ref{Sect3}, $\psi(\lambda)$ is of the above form with $h=N$. So
$\psi^\ast(\lambda)$ is defined in each case. By \cite{Saito94} for the
exceptional unimodal singularities and \cite{Ebeling96} for the extended
duality we have

\begin{theorem} \label{thm:Saito}
If $(X,0)$ is a singularity of class a), b), c), or d) with
associated polynomial $\psi(\lambda)$, then $\psi^\ast(\lambda)$ is the
corresponding polynomial of the dual singularity.
\end{theorem}

This is proved in \cite{Saito94} and \cite{Ebeling96} by a case-by-case
verification. Here we derive this result
from the duality of Newton polyhedra considered in Sect.~\ref{Sect2}.

\addvspace{3mm}

\noindent {\em Proof of Theorem~\ref{thm:Saito}. }
Let $\Gamma$ be the Newton diagram of $(X,0)$ and $\nabla$ be the dual
polygon.

First assume that $(X,0)$ is of type a) or b). Then the polygon
$\nabla$ has the following vertices:
$$\tilde{v}_1=(w_1+1,w_2,w_3), \tilde{v}_2=(w_1,w_2+1,w_3),
\tilde{v}_3=(w_1,w_2,w_3+1).$$
We have
$$V(\nabla) = \left|
\begin{array}{ccc} w_1+1 & w_1 & w_1 \\
                   w_2 & w_2+1 & w_2 \\
                   w_3 & w_3 & w_3+1
\end{array} \right| = w_1 + w_2 +w_3 + 1 = N,$$
$$\tilde{v}_1 \times \tilde{v}_2 = \left( \begin{array}{c}
- w_3 \\ -w_3 \\ N-w_3
\end{array} \right),
\quad \tilde{v}_1 \times \tilde{v}_3 = \left( \begin{array}{c}
 w_2 \\ -N+w_2 \\ w_2
\end{array} \right),
\quad \tilde{v}_2 \times \tilde{v}_3 = \left( \begin{array}{c}
N- w_1 \\ -w_1 \\ -w_1
\end{array} \right).$$

The exponent of $(\lambda-1)$ in $\psi(\lambda)=\phi(\lambda)$ is $\chi_1=-1$
by Sect.~\ref{Sect3}. Since
$V(\nabla)= N$ the exponent of $(\lambda^N-1)$ in $\psi^\ast(\lambda)$
is equal to $-\chi_1=1$.

Let $s$ be a special 0-dimensional face of $\Gamma$. Because $s$ is special,
it lies on one of the coordinate axes.  Without loss of generality we assume
that $s$ lies on the $x$-axis. (The other cases are treated similarly.) Then
$s=(\frac{N}{w_1},0,0)$. This contributes a factor $(\lambda^{N/w_1} -1)$ to
$\phi(\lambda)$. In particular, $w_1 | N$, so $N=aw_1$ for some integer $a
\geq 2$. Then $\tilde{v}_2 \times \tilde{v}_3 = w_1 (a-1, -1, -1)^t$. If the
1-simplex $\mbox{pos}\{\tilde{v}_2,\tilde{v}_3\}$ is special, i.e.,
lies in a coordinate hyperplane (with respect
to the basis given by the rows of $A$), then we have a contribution
$(\lambda^{w_1} -1)^{-1}$ to $\phi^\ast(\lambda)$.

Now let $s$ be a special 1-dimensional face of $\Gamma$. We assume without
loss of generality that $s$ lies in the $x,y$-plane. Let $v = (x,y,0)$ and
$v'=(x',y',0)$ be the vertices of $s$. Let $d=\gcd(w_1,w_2)$ and $\delta=
\gcd(x-x',y-y')$. Then $m(s)=\frac{N}{d}$ and
$V(s) = xy'-x'y = \frac{N}{d} \delta= m(s)\delta$. This contributes a factor
$(\lambda^{N/d}-1)^{-\delta}$ to $\phi(\lambda)$. On the other hand
$\tilde{v}_3 = (w_1, w_2, N-w_1-w_2) = d
(\frac{w_1}{d},\frac{w_2}{d},\frac{N}{d}-\frac{w_1}{d}-\frac{w_2}{d})$. If
$\frac{1}{d}\tilde{v}_3$ is a row of $A$, then we have a contribution
$(\lambda^d-1)$ to $\phi^\ast(\lambda)$.

Now consider a singularity $(X,0)$ of class c). The polygon
$\nabla$ has the following vertices:
$$\tilde{v}_1=(w_1+2,w_2+1,w_3), \tilde{v}_2=(w_1,w_2+1,w_3),
\tilde{v}_3=(w_1,w_2,w_3+1).$$
Let $N_1$ be the degree of the polynomial $g(x,y,z,w) = xw-y^2$. Then $N_1 =
w_1+w_4 = 2w_2$. Moreover, we have
$$w_1+w_2+w_3+w_4+1=N_1+N, \quad w_3=\frac{N}{2}.$$
This implies that $N=w_2+w_3+1$. Therefore we have
$$V(\nabla) = \left|
\begin{array}{ccc} w_1+2 & w_1 & w_1 \\
                   w_2+1 & w_2+1 & w_2 \\
                   w_3 & w_3 & w_3+1
\end{array} \right| = 2(w_2 +w_3 + 1) = 2 N,$$
$$\tilde{v}_1 \times \tilde{v}_2 = 2 \left( \begin{array}{c}
0 \\ -w_3 \\ N-w_3
\end{array} \right),
\quad \tilde{v}_1 \times \tilde{v}_3 = \left( \begin{array}{c}
 N \\ w_4 -2N \\ w_4
\end{array} \right),
\quad \tilde{v}_2 \times \tilde{v}_3 = \left( \begin{array}{c}
N \\ -w_1 \\ -w_1
\end{array} \right).$$

The exponent of $(\lambda-1)$ in $\psi(\lambda)$ is $\chi_1=-2$. Since
$V(\nabla)= 2N$ the exponent of $(\lambda^N-1)$ in $\psi^\ast(\lambda)$
is equal to $-\chi_1=2$.

Let $s$ be a special $0$-dimensional face of $\Gamma$. Then $s$ lies on one
of the coordinate axes $\RR e_1$, $\RR (-e_1+2e_2)$, or $\RR e_3$. Then we
have $s=\frac{N}{w_1}e_1$, $s=\frac{N}{w_4}(-e_1+2e_2)$, or
$s=\frac{N}{w_3}e_3=2e_3$ respectively. This contributes a factor
$(\lambda^{N/w_i}-1)$ to $\psi(\lambda)$ where
$i=1,4,3$ respectively. For $i=3$ this factor is equal to $(\lambda^2-1)$.
Theorem~\ref{thm:EVc} implies that
$$\phi(\lambda) = (\lambda^2-1) \phi_{\Gamma(p)}(\lambda)$$
where $p(x,y,z)=f(x,y,z,x^{-1}y^2)$.
From the above formulae for the mutual
vector products of the vectors $\tilde{v}_1$, $\tilde{v}_2$, $\tilde{v}_3$ we
see that we have a possible contribution
$(\lambda^{w_1}-1)^{-1}$,
$(\lambda^{w_4}-1)^{-1}$, or $(\lambda^{w_3}-1)^{-2}$ respectively to
$\psi^\ast(\lambda)$.

Now let $s$ be a special 1-dimensional face of $\Gamma$. Let
$v=x_1e_1+x_2(-e_1+2e_2)+x_3e_3$ and $v'=x'_1e_1+x'_2(-e_1+2e_2)+x'_3e_3$ be
the vertices of $s$. Since $s$ lies in a coordinate plane, $x_i=x'_i=0$ for
some $i=1,2,3$. As above, we see that we have a contribution
$(\lambda^{N/d_i}-1)^{-\delta_i}$ to $\psi(\lambda)$, where
\begin{eqnarray*}
d_1= \gcd(w_3,w_4), & \quad & \delta_1=\gcd(x_2-x_3,x'_2-x'_3),\\
d_2= \gcd(w_1,w_3), & \quad & \delta_2=\gcd(x_1-x_3,x'_1-x'_3),\\
d_3= \gcd(w_1,w_4), & \quad & \delta_3=\gcd(x_1-x_2,x'_1-x'_2).
\end{eqnarray*}
On the other hand
$$\tilde{v}_1 = \left( \begin{array}{c} N-w_4 \\ N-w_3 \\ w_3 \end{array}
\right), \quad \tilde{v}_2 = \left( \begin{array}{c} w_1 \\ N-w_3 \\ w_3
\end{array} \right), \quad \tilde{v}_3 = \left( \begin{array}{c} w_1 \\
\frac{w_1+w_4}{2} \\ N-\frac{w_1+w_4}{2} \end{array}
\right).$$
This shows that if $\frac{1}{d_i}\tilde{v}_i$ is a row of $A$ and 2 does
not divide $d_3$ then there is a factor $(\lambda^{d_i}-1)$ in
$\psi^\ast(\lambda)$.

Finally let $(X,0)$ be of type d). The polygon
$\nabla$ has the following vertices:
$$\tilde{v}_1=(w_1+1,w_2,w_3+1), \tilde{v}_2=(w_1+1,w_2+1,w_3),
\tilde{v}_3=(w_1,w_2,w_3+1),$$
$$\tilde{v}_4=(w_1,w_2+1,w_3).$$
Let $N_1$ be the degree of the polynomial $g(x,y,z,w) = xw-yz$. Then $N_1 =
w_1+w_4 = w_2+w_3$. Moreover, we have
$$w_1+w_2+w_3+w_4+1=N_1+N.$$
This again implies that $N=w_2+w_3+1$. Therefore we have
\begin{eqnarray*}
V(\nabla)  & = & \left|
\begin{array}{ccc} w_1+1 & w_1+1 & w_1 \\
                   w_2 & w_2+1 & w_2 \\
                   w_3+1 & w_3 & w_3+1
\end{array} \right| + \left|
\begin{array}{ccc} w_1+1 & w_1 & w_1 \\
                   w_2+1 & w_2+1 & w_2 \\
                   w_3 & w_3 & w_3+1
\end{array} \right| \\
& = & 2(w_2+w_3+1) \\
& = & 2N.
\end{eqnarray*}

The exponent of $(\lambda-1)$ in $\psi(\lambda)$ is $\chi_1=-2$. Since
$V(\nabla)= 2N$, the exponent of $(\lambda^N-1)$ in $\psi^\ast(\lambda)$
is equal to $-\chi_1=2$.

Let $s$ be a special $0$-dimensional face of $\Gamma$. By inspection of
the equations of Table~\ref{Table1} we see that $s$ either lies on
$\RR e_1$ or on $\RR e_3$. Then we have $s=\frac{N}{w_1}e_1$ or
$s=\frac{N}{w_3}e_3$ respectively. Then $\psi(\lambda)$ has a factor
$(\lambda^{N/w_1}-1)$ or $(\lambda^{N/w_3}-1)$ respectively. We have
$$\tilde{v}_3 \times \tilde{v}_4 = \left( \begin{array}{c}
-N \\ w_1 \\ w_1 \end{array} \right),
\quad \tilde{v}_2 \times \tilde{v}_4 = \left( \begin{array}{c}
 0 \\ -w_3 \\ N-w_3
\end{array} \right).$$
If the 1-simplices $\mbox{pos}\{\tilde{v}_3,\tilde{v}_4\}$ or
$\mbox{pos}\{\tilde{v}_2,\tilde{v}_4\}$ are special then we have a contribution
$(\lambda^{w_1} -1)^{-1}$ or $(\lambda^{w_3} -1)^{-1}$ respectively to
$\psi^\ast(\lambda)$.

Now let $s$ be a special 1-dimensional face of $\Gamma$. By inspection of the
equations of \cite{Ebeling96} one sees that $s$ can only lie in the
$x,y$-plane or $x,z$-plane. In the first case we have a factor
$(\lambda^{N/d}-1)^{-\delta}$ of $\psi(\lambda)$ where $d = \gcd(w_1,w_2)$
and $\delta$ is as above. On the other hand, $\tilde{v}_3 = (w_1, w_2,
N-w_2)$, so if $\frac{1}{d}\tilde{v}_3$ is a row of $A$, then we have a
factor $(\lambda^d-1)$ in $\psi^\ast(\lambda)$. The other case is treated
similarly using $\tilde{v}_4 = (w_1, N-w_3, w_3)$.

This finishes the proof of Theorem~\ref{thm:Saito}.

\addvspace{10mm}

{\sc Institut f\"{u}r Mathematik, Universit\"{a}t Hannover, Postfach 6009,
D-30060 Hannover, Germany}

{\em E-mail address}: ebeling@math.uni-hannover.de
\end{document}